\documentclass[12pt,reqno]{amsart}
\usepackage{amsmath,amssymb}
\usepackage[dvips]{graphicx,color}
\usepackage{latexsym,amssymb}
\usepackage{mathrsfs}
\usepackage[abbrev]{amsrefs}

\newtheorem{thm}{Theorem}[section]
\newtheorem{cor}[thm]{Corollaly}
\newtheorem{prop}[thm]{Proposition}
\newtheorem{lem}[thm]{Lemma}

 \newenvironment{pf}
    {{\noindent \bf Proof. }}{\hfill $\Box$}
\numberwithin{equation}{section}
\numberwithin{thm}{section}
\begin{document}

\begin{center}
 \textbf{\Large{The Derivative Structure for a Quadratic Nonlinearity 
 and Uniqueness for SQG}}
\end{center}

\vskip5mm 

\centerline{Tsukasa Iwabuchi$^*$ }
\centerline{Mathematical Institute, Tohoku University}
\centerline{Sendai 980-8578 Japan}
\footnote[0]{E-mail: $^*$t-iwabuchi@tohoku.ac.jp}

\begin{center}
 \begin{minipage}{120mm}
\small \textbf{Abstract.} 
We study the two-dimensional surface quasi-geostrophic equation on a bounded 
domain with a smooth boundary. Motivated by the three-dimensional incompressible 
Navier-Stokes equations and previous results in the entire space $\mathbb R^2$, 
we demonstrate that the uniqueness of the mild solution holds in $L^2$. 
For the proof, 
we provide a method for handling fractional Laplacians in nonlinear problems, 
and  
develop an approach to derive second-order derivatives 
for the nonlinear term involving fractional derivatives 
of the Dirichlet Laplacian.
 \end{minipage}
\end{center}

\section{Introduction}

We study bilinear estimates related to the Leibniz rule for certain quadratic 
nonlinearities. The simplest example is given by
\[
(\nabla f) g + f \nabla g = \nabla (fg).
\]
Note that each term on the left-hand side is not a derivative of a single function, while the right-hand side represents the derivative of the first order.

Another example comes from the Keller-Segel equation of the parabolic-elliptic type:
\[
\partial_t u - \Delta u = \nabla \cdot (u \nabla \psi), \qquad -\Delta \psi = u.
\]
We can express this as
\[
\begin{split}
\nabla \cdot (u \nabla \psi) &= \nabla \cdot \left( u \nabla (-\Delta)^{-1} u \right) \\
&= - \sum_{j,k} \nabla_j \nabla_k \cdot \left\{ \left( \nabla_k (-\Delta)^{-1} u \right) \left( \nabla_j (-\Delta)^{-1} u \right) \right\} + \Delta \left| \nabla (-\Delta)^{-1} u \right|^2.
\end{split}
\]
Thus, the nonlinear term can be written in terms of second-order derivatives.

Another example is the two-dimensional surface quasi-geostrophic (SQG) equation:
\begin{equation}\label{0830-1}
\partial_t \theta - \Delta \theta + (u \cdot \nabla) \theta = 0, \qquad u = \nabla^\perp \Lambda^{-1} \theta = (-\partial_{x_2}, \partial_{x_1}) (-\Delta)^{-\frac{1}{2}} \theta,
\end{equation}
where $\Lambda = (-\Delta)^{\frac{1}{2}}$. The nonlinear term can also be expressed using second-order derivatives. Specifically,
\[
(u \cdot \nabla) \theta = \nabla \cdot \left( \left(\nabla^\perp \Lambda^{-1} \theta \right) \theta \right),
\]
and
\begin{equation}\label{0830-3}
\begin{split}
& (\nabla^\perp \Lambda^{-1} f) g + (\nabla^\perp \Lambda^{-1} g) f \\
=& {\mathcal F}^{-1} \left[ \int_{{\mathbb R}^2} i \left( \frac{(\xi - \eta)^\perp}{\lvert \xi - \eta \rvert} + \frac{\eta^\perp}{\lvert \eta \rvert} \right) \hat{f}(\xi - \eta) \hat{g}(\eta) \mathrm{d}\eta \right] \\
=& i {\mathcal F}^{-1} \left[ \int_{{\mathbb R}^2} \frac{\xi^\perp}{\lvert \xi - \eta \rvert} \hat{f}(\xi - \eta) \hat{g}(\eta) \mathrm{d}\eta \right] \\
&+ i {\mathcal F}^{-1} \left[ \int_{{\mathbb R}^2} \frac{\xi \cdot (\xi - 2\eta)}{\lvert \xi - \eta \rvert + \lvert \eta \rvert} \cdot \frac{1}{\lvert \xi - \eta \rvert} \hat{f}(\xi - \eta) \cdot \frac{\eta^\perp}{\lvert \eta \rvert} \hat{g}(\eta) \mathrm{d}\eta \right] \\
=& \nabla^\perp \left\{ (\Lambda^{-1} f) g \right\} + \nabla \cdot m(D_1, D_2)\left\{ (\Lambda^{-1} f) (\nabla^\perp \Lambda^{-1} g) \right\},
\end{split}
\end{equation}
where $\displaystyle 
(m(D_1, D_2)(fg) := {\mathcal F}^{-1} \left[ \int_{{\mathbb R}^2} \frac{\xi - 2\eta}{\lvert \xi - \eta \rvert + \lvert \eta \rvert} \hat{f}(\xi - \eta) \hat{g}(\eta) \mathrm{d}\eta \right]$. We then see that the nonlinear term \((u \cdot \nabla) \theta\) can be written in terms of second-order derivatives.

We refer to \cites{Iw-2011,IwNa-2013,OhIw-2022} for the ideas behind these derivatives.
Thanks to this structure, we can prove several bilinear estimates and demonstrate the existence of solutions using the fixed-point theorem 
and large time behavior.

An advantage of the derivative forms discussed in this paper is that they allow for valid estimates of norms in Besov spaces $\dot{B}^0_{1,q}$. Without such derivative forms, we can find functions $f, g \in \mathcal{S}(\mathbb{R}^d)$ such that
\[
\| fg \|_{\dot{B}^0_{1,q}} = \infty, \quad \text{if } q < \infty.
\]
The Gaussian is one of the typical examples of the above. We refer to the book~\cite{RuSi_1996} (see subsection~4.8.3) for bilinear estimates in Besov spaces with the regularity index $s = 0$.

On the other hand, the derivative forms allow us to obtain product estimates thanks to the following inequality. If $s > 0$, then
\[
\| fg \|_{\dot{B}^s_{1,q}} \leq C \| f \|_{\dot{B}^s_{p_1,q}} \| g \|_{L^{p_2}} + C \| f \|_{L^{p_3}} \| g \|_{\dot{B}^s_{p_4,q}},
\]
with $1 = 1/p_1 + 1/p_2 = 1/p_3 + 1/p_4$. 

Here, we focus on the quadratic nonlinearity appearing in the surface quasi-geostrophic equation~\eqref{0830-1}. In the Euclidean space $\mathbb{R}^2$, by using the Fourier transformation, it is possible to prove that
\[
\| (\nabla (-\Delta)^{-\frac{1}{2}} \theta ) \theta \|_{\dot{B}^0_{1,q}} \leq C \| \theta \|_{\dot{B}^0_{p_1,q}} \| \theta \|_{L^{p_2}}.
\]
To study the product of two functions, we consider
\begin{equation}\label{0830-2}
(\nabla \Lambda^{-1} f ) g + f \nabla \Lambda^{-1} g
\end{equation}
instead of $(\nabla (-\Delta)^{-\frac{1}{2}} \theta ) \theta$ 
as in \cites{OhIw-2022,IwUe-2024}.

In this paper, we study the inequality for such quadratic terms. To this end, 
we need to introduce the fractional Laplacian and consider the Dirichlet 
Laplacian on $L^2$ in a smooth bounded domain. Throughout this paper, 
let $\Omega$ be a smooth bounded domain in $\mathbb{R}^2$. We denote 
the Dirichlet Laplacian on $L^2$ by $-\Delta_D$. By the arguments 
in \cites{FaIw-2024,IMT-2019}, we regard $-\Delta_D$ as operators 
on distribution spaces and Besov spaces, whose definitions are given 
in Section~\ref{sec:2}.

We will define the Dirichlet Laplacian $-\Delta_D$ on distribution spaces 
denoted by $\mathcal{Z}'$ and introduce the Besov spaces associated 
with the Dirichlet Laplacian. Let $\phi$ satisfy
\[
\phi \in C_0^\infty (\mathbb{R}), 
\quad 
\text{supp } \phi \subset [2^{-1}, 2], 
\quad 0 \leq \phi \leq 1, 
\quad 
\sum_{j \in \mathbb{Z}} \phi (2^{-j}\lambda) = 1, 
\quad 
\lambda > 0.
\]
We introduce $\{ \phi_j \}_{j \in \mathbb{Z}}$ such that for every $j \in \mathbb{Z}$
\[
\phi_j(\lambda) := \phi(2^{-j} \lambda), \quad \lambda \in \mathbb{R},
\]
which forms a partition of unity by dyadic numbers on the half-line 
$(0, \infty)$. Formally, for $s \in \mathbb{R}$ and $1 \leq p, q \leq \infty$, 
the norm of the Besov space $\dot{B}^s_{p,q}$ is given by
\[
\| u \|_{\dot{B}^s_{p,q}} := 
\left\| \left\{ 2^{js} \| \phi_j(\sqrt{-\Delta _D}) u \|_{L^p} 
\right\}_{j \in \mathbb{Z}} 
\right\|_{\ell^q (\mathbb{Z})},
\]
where $\phi_j (\sqrt{-\Delta _D})$ denotes the operator associated with $\phi_j$ 
and $\sqrt{-\Delta _D}$. We also denote the square root of the Dirichlet Laplacian 
by
\[
\Lambda_D = \sqrt{-\Delta_D}.
\]
We will provide the precise definition of $\dot{B}^s_{p,q}$ 
in Section~\ref{sec:2}.

\vskip3mm 

This paper is organized as follows. In subsection~\ref{subsec:1.1}, 
we present the first theorem on the bilinear estimate for the product \eqref{0830-2} involving the square root of the Dirichlet Laplacian. 
In subsection~\ref{subsec:1.2}, we discuss the second theorem on 
the uniqueness problem for the two-dimensional surface quasi-geostrophic 
equation in $C([0,T], L^2)$ on a smooth domain. We will develop 
the argument for the product \eqref{0830-2} to derive derivative forms 
for the nonlinear term. In section~\ref{sec:2}, we introduce the definition 
of the Besov spaces associated with the Dirichlet Laplacian and prepare 
several lemmas and propositions for the proofs of the theorems. 
In sections~\ref{sec:3} and~\ref{sec:4}, we prove the first and 
second theorems, respectively.

\subsection{Bilinear Estimate}\label{subsec:1.1}

The following is the first theorem of the present paper regarding the bilinear estimate of the product as in \eqref{0830-2}, where $(-\Delta)^{\frac{1}{2}}$ is replaced by the square root of the Dirichlet Laplacian $\Lambda_D = \sqrt{-\Delta_D}$. Let $\dot{B}^s_{p,q}$ denote the Besov spaces associated with the Dirichlet Laplacian.

\begin{thm}\label{Thm:1}
Let $\Omega$ be a bounded domain in $\mathbb{R}^2$ with a smooth boundary. 
Let $-1 < s < 2$, $1 \leq p, p_1, p_4, q \leq \infty$, 
$1 < p_2, p_3 < \infty$, and  
\[
\dfrac{1}{p} = \dfrac{1}{p_1} + \dfrac{1}{p_2} = \dfrac{1}{p_3} + \dfrac{1}{p_4}.
\]
Then there exists a constant $C > 0$ such that for every $f \in \dot{B}^s_{p_1,q} \cap L^{p_3}$ and $g \in L^{p_2} \cap \dot{B}^s_{p_4,q}$, the function $(\nabla \Lambda_D^{-1} f) g + f \nabla \Lambda_D^{-1} g$ belongs to $\dot{B}^s_{p,q}$ and 
\[
\| (\nabla \Lambda_D^{-1} f) g + f \nabla \Lambda_D^{-1} g \|_{\dot{B}^s_{p,q}} \leq C \left( \| f \|_{\dot{B}^s_{p_1,q}} \| g \|_{L^{p_2}} + \| f \|_{L^{p_3}} \| g \|_{\dot{B}^s_{p_4,q}} \right).
\]
\end{thm}

\noindent 
{\bf Remark. }  
(1) The regularity restriction $s > -1$ is essential 
because $(\nabla \Lambda_D^{-1} f) g + f \nabla \Lambda_D^{-1} g$ 
is expressed in terms of the derivative form of the first order. 
The condition $s < 2$ could be improved to $s < 2 + \frac{1}{p}$, 
as discussed in~\cite{Iw-2023}, but this is not the focus of the present paper.

\vskip2mm   

\noindent 
(2) We mention the case with a more general power of the 
Laplacian. Let us consider the Euclidean space $\mathbb{R}^2$ and 
\[
\Big(\nabla \Phi(D) f \Big) g + f \Big(\nabla \Phi(D) g \Big),
\]
where $\Phi (\lambda)= \lambda^{\alpha}$ ($\alpha \in \mathbb{R}$), and 
$\Phi(D)f = \mathcal{F}^{-1} \Phi(|\xi|) \mathcal{F} f$. By the Fourier transformation, 
\[
\begin{split}
\mathcal{F}\Big[ 
\Big(\nabla \Phi(D) f \Big) g + f \Big(\nabla \Phi(D) g \Big) \Big]
=
& i\int_{\mathbb{R}^d}
  \Big\{ (\xi - \eta) \Phi(\xi - \eta) + \eta \Phi(\eta) \Big\} 
   f(\xi - \eta) g(\eta) \, {\rm d}\eta,
\end{split}
\]
and by the mean value theorem,
\[
\begin{split}
  (\xi - \eta) \Phi(\xi - \eta) + \eta \Phi(\eta) 
=& \xi \Phi(\xi - \eta) + \eta (\Phi(\eta) - \Phi(\eta -\xi)) \\
=& \xi \Phi(\xi - \eta) + \eta \xi \cdot \nabla \Phi (\eta - \theta \xi) 
\quad \text{for some } \theta \in (0,1),
\end{split}
\]
which also yields a derivative form due to the multiplication by $\xi$. 
Therefore, we can expect the inequality with a more general fractional Laplacian.

\vskip3mm

We give a few comments on the proof of Theorem~\ref{Thm:1}. 
The crucial point 
to derive the derivative of the 1st order for 
$ (\nabla \Lambda_D^{-1} f ) g + f \nabla \Lambda_D^{-1} g$ 
requires a more complicated argument 
compared with \eqref{0830-3} in the $\mathbb{R}^2$ case. 
To this end, we use the formula of the square root of the Dirichlet Laplacian 
by the resolvent,
\begin{equation}\label{0904-4}
\Lambda _D = c_0 \int_0 ^\infty \mu^{-\frac{3}{2}} 
 \Big( 1 - (1- \mu \Delta _D)^{-1} \Big) \, {\rm d}\mu,
\end{equation}
which allows us to have a derivative form 
(see~\eqref{0830-7} and below). 
Once we obtain the form, it is possible to estimate the Besov norm 
using bilinear estimates in function spaces with positive regularity, 
along with second-order derivative estimates.

We also note from our proof below that, instead of the Bony paraproduct 
formula in~\cite{Bo-1981},  
we use a simple decomposition for the product:
\[
fg = \Big( \sum_{k\geq l} + \sum_{k < l}\Big) 
 \Big( \phi_j (\Lambda_D) f\Big)  \Big( \phi_j (\Lambda_D) g\Big).
\]
This approach works for our proof.

\subsection{Uniqueness for SQG} \label{subsec:1.2}

We study the two-dimensional surface quasi-geostrophic equation. 
Let $\Omega$ be a bounded domain with a smooth boundary. 

\begin{equation}\label{eq:SQG}
 \left\{
   \begin{aligned}
     &\partial_t \theta + (u \cdot \nabla) \theta - \Delta_D \theta = 0,
           \quad &&\text{in }  (0,T) \times \Omega, \\
     &u = \nabla^{\perp} \Lambda_D^{-1} \theta,
           \quad &&\text{in }  (0,T) \times \Omega, \\
     &\theta = 0,
           \quad &&\text{on }  (0,T) \times \partial \Omega, \\
     &\theta(0,\cdot) = \theta_0(\cdot),
           \quad &&\text{in } \Omega,
   \end{aligned}
  \right.
\end{equation}
where $-\Delta_D$ is the Dirichlet Laplacian, 
$\nabla^{\perp} = (- \partial_{x_2}, \partial_{x_1})$, and $\Lambda_D = (-\Delta_D)^{\frac{1}{2}}$. 
The equations are derived from general quasi-geostrophic equations in the special case of constant potential vorticity and buoyancy frequency (see~\cite{Pedlosky}). 
It is known that they describe atmospheric motion and are useful in weather forecasting. 
Here, $\theta$ represents potential temperature, and $u$ denotes the atmospheric velocity vector.

\vskip2mm

We introduce the definition of mild solutions and state the second theorem.

\vskip3mm

\noindent \textbf{Definition.} Let $T > 0$ and $\theta_0 \in L^2$.  
If $\theta : [0,T] \times \Omega \to \mathbb{R}$ 
satisfies
\[
\begin{cases}
\theta \in C([0,T] , L^2), \\
\displaystyle 
\int_{\Omega} \theta (t) g \, \mathrm{d}x
= \int_{\Omega} (e^{t\Delta_D} \theta_0) g \, \mathrm{d}x 
+ \int_0^t \int_{\Omega} u \theta \cdot \nabla e^{(t-\tau)\Delta_D} g 
  \, \mathrm{d}x \, \mathrm{d}\tau, 
\quad \text{for all } t \in [0,T], g \in \mathcal{Z},
\end{cases}
\]
then we call $\theta$ a mild solution of \eqref{eq:SQG}. 
Here, the space $\mathcal{Z}$ is a Fr\'echet space defined  
by \eqref{0902-1} in section~\ref{sec:2}, and it consists of 
$g \in L^2$ such that $(-\Delta_D)^M g \in L^1$ for all 
natural numbers $M$.

\vskip3mm 

It is a classical result that for every $\theta_0 \in L^2$, 
there exists a unique global mild solution $\theta$ such that 
$\theta \in C([0,\infty);L^2) \cap L^2(0,\infty; \dot H^1)$. 
Our aim is to establish uniqueness without the condition $L^2(0,\infty; \dot H^1)$.

\begin{thm}\label{Thm:2}
Let $T > 0$. Suppose that $\theta, \tilde{\theta} \in C([0,T];L^2)$ are  
mild solutions of \eqref{eq:SQG} such that 
$\theta (0) = \widetilde \theta (0) \text{ in } L^2$. 
Then $\theta (t) = \widetilde \theta (t) \text{ in } L^2$  
for all $t \in [0,T]$. 
\end{thm}

Let us recall several known results about the uniqueness of \eqref{eq:SQG}. 
As mentioned before, the uniquness for 
$u \in L^{\infty}(0,T ; L^2) \cap L^2(0,T; \dot{H}^1)$ is well-known,   
where the domain is either $\mathbb R^2$ 
or a bounded domain with a smooth boundary. 
The uniqueness in $C([0,T], L^2)$ for $\Omega = \mathbb R^2$ 
was proved in the paper~\cite{IwUe-2024}. 
We also refer to several results~\cites{CoWu-1999,Fe-2011}  
for the fractional Laplacian 
$(-\Delta)^{\frac{\alpha}{2}}$, $1 < \alpha \leq 2$, on $\mathbb R^2$. 
The scale critical spaces of the two-dimensional surface quasi-geostrophic 
equation are the same as those of the two-dimensional incompressible 
Navier-Stokes equations. 
We also refer to the paper \cite{ChLu-2023} by Cheskidov-Luo, 
who proved the non-uniqueness in $C([0, T] ; L^p)$ with $p<2$ 
for the two-dimensional incompressible Navier-Stokes equations. 
The case when $p = 2$ does not seem to be resolved.

Regarding the three-dimensional incompressible Navier-Stokes equations, 
Kato~\cite{Ka-1984} showed the uniqueness in 
the class $C([0,T), L^3) \cap C((0,T), L^r)$, with $3 < r \leq \infty$, under the 
convergence condition 
$t^{\frac{3}{2}(\frac{1}{3}-\frac{1}{r})} \| u(t) \| \to 0 \, (t\to 0)$. 
Brezis~\cite{Br-1994} removed the condition by showing that 
every mild solution in the class satisfies the convergence. 
Meyer~\cite{May-1997}, Furioli, Lemari\'{e}-Rieusset, and Terraneo~\cite{FLT-1997} and Monniaux~\cite{Monn-1999} demonstrated that the uniqueness of the mild solution $u \in C ([0, T ]; L^3 )$ holds. Lions and Masmoudi~\cite{LiMa-2001} proved the uniqueness in $L^3$ by solving a dual problem. Furioli, Lemari\'{e}-Rieusset, and Terraneo~\cite{FLT-1997} used Besov spaces, and the proof of the main theorem in this paper is inspired by their result. 
We also refer to an important result 
on non-uniqueness by Buckmaster and Vicol~\cite{BuVi-2019}, 
who proved the non-uniqueness in the class $L^\infty (0,T ; L^2)$.

Furthermore, we have the following corollary since there is a global solution 
$\theta \in C([0,\infty) ; L^2) \cap L^2(0,\infty; \dot{H}^1)$ 
for all $\theta_0 \in L^2$.

\begin{cor}
Consider the equations \eqref{eq:SQG}  
in $(0, \infty) \times \Omega$. Then, 
for every $\theta_0 \in L^2$, there is a unique global mild solution 
$\theta$ such that $\theta \in C([0, \infty) ; L^2)$.
\end{cor}

We provide a few comments on the proof. We establish the argument, 
as done in \cite{IwUe-2024}, involving the fractional Dirichlet Laplacian.
The crucial point in the proof of Theorem~\ref{Thm:2} is 
to derive the second-order derivative for the nonlinear term 
$(u \cdot \nabla) \theta$ involving the fractional Dirichlet Laplacian. 
However, the fractional Dirichlet Laplacian is defined via the 
spectral theorem, which we cannot handle similarly to the Fourier 
transform. 
To address the issue, we use the formula \eqref{0904-4},  
which allows us to have a derivative form. 
We also need to give a definition for the derivative in the nonlinear term, 
since the gradient itself changes the boundary condition in general.

\vskip3mm

\section{Preliminaries} \label{sec:2}

We introduce test function spaces, distribution spaces, 
Besov spaces associated with the Dirichlet Laplacian, 
together with the duality property, as in the paper~\cite{IMT-2019}. 
We also provide the definition of derivatives for functions belonging to $L^1$.

\subsection{Besov Spaces and Basic Properties}
\quad 

\vskip3mm 

\noindent 
{\bf Definition. } 
We define $-\Delta _D$ by the Dirichlet Laplacian on $L^2$, i.e., 
\[
\begin{cases}
D(-\Delta _D) = \{ f \in H^1_0 (\Omega) \mid  \Delta f \in L^2 \},
\\
-\Delta _D f = - (\partial_{x_1}^2 + \partial_{x_2}^2) f , 
\quad f \in D(-\Delta _D), 
\end{cases}
\]
and $\Lambda _D := \sqrt{-\Delta _D}$.  
We also note that the above definition is concerned with $L^2$, but 
we use the same notation $-\Delta _D$ for several spaces, 
such as Besov spaces and distribution spaces, appearing in this 
paper.

\vskip3mm

The following lemma is about the boundedness of the spectral multiplier, 
which is essential for extending the operator from $L^2$ to $L^p$. 
The proof idea is based on the work by Jensen and Nakamura~\cite{JeNa-1995}, 
who investigated this in the Euclidean space. 
We also refer to the book by Ouhabaz~\cite{Ou_2005} and 
the paper~\cite{IMT-2018} for general domains.

\begin{lem}{\rm (}\cite{IMT-2018}{\rm )} 
\label{lem:0809-1}
Let $1 \leq p \leq \infty$, $\phi \in C_0^\infty (\mathbb R)$, 
and $\phi_j(\lambda ) = \phi (2^{-j}\lambda)$ for $j \in \mathbb Z$ 
and $\lambda \in \mathbb R$. 
Then we have 
\[
\sup_{j \in \mathbb Z} \| \phi_j(\Lambda_D) \|_{L^p \to L^p} < \infty. 
\]
Furthermore, if $1 \leq r \leq p \leq \infty$, then 
\[
\sup_{j \in \mathbb Z} 
2^{-\alpha j -n\left(\frac{1}{r}-\frac{1}{p}\right)j}\| \Lambda^\alpha \phi_j(\Lambda_D) \|_{L^r \to L^p} 
< \infty. 
\]
\end{lem}

\vskip3mm

Following the procedure in the paper~\cite{IMT-2019} (see also \cite{FaIw-2024}), 
we introduce test function spaces $\mathcal X$ and $\mathcal Z$, 
which are Fr\'echet spaces,  
and spaces of distributions as the topological duals of these. 

\vskip3mm 

\noindent 
{\bf Definition. } (Test function spaces and distributions)
\begin{enumerate}
\item[(i)] (Non-homogeneous type) 
We define a test function space $\mathcal X$ by 
\[
\mathcal X
:= \{ u \in L^1 \cap L^2 \mid p_M(u) < \infty \text{ for all } M \in \mathbb N\}, 
\]
where 
\[
p_M(u) := \| u \|_{L^1} + \sup_{j \in \mathbb N} 2^{Mj} \| \phi_j(\Lambda_D) u \|_{L^1}. 
\]

\item[(ii)] (Homogeneous type) 
We define a test function space $\mathcal Z$ by 
\begin{equation}\label{0902-1}
\mathcal Z 
:= \{ u \in \mathcal X \mid q_M(u) < \infty 
 \text{ for all } M \in \mathbb N\}, 
\end{equation}
where 
\[
q_M(u) := p_M(u) + \sup_{j \leq 0} 2^{M|j|} \| \phi_j(\Lambda_D) u \|_{L^1}. 
\]

\item[(iii)] 
$\mathcal X', \mathcal Z'$ are defined by the topological duals of 
$\mathcal X, \mathcal Z$, respectively. 

\end{enumerate}

\vskip3mm

\noindent 
{\bf Remark. } 
$\mathcal X, \mathcal X'$ are function spaces for the non-homogeneous spaces, 
and $\mathcal Z, \mathcal Z'$ are for the homogeneous spaces 
(see~\cites{FaIw-2024,IMT-2019}). 
Here the domain we consider is a bounded domain, and the infimum of the 
spectrum of $-\Delta _D$ is strictly positive. 
Therefore 
\[
\mathcal X \simeq \mathcal Z \quad \text{and}\quad 
\mathcal X' \simeq \mathcal Z'. 
 \]
In what follows, we choose $\mathcal Z, \mathcal Z'$, 
to define Besov spaces of the homogeneous type. 

\vskip3mm 

We give the following definition of how to understand $\nabla f $ 
as an element of our distribution space $\mathcal Z'$. 

\vskip3mm 

\noindent 
{\bf Definition. } 
Let $f \in L^1$. We define $\nabla f$ as an element of $\mathcal Z'$ 
by 
\begin{equation}\label{0904-8}
{}_{\mathcal Z'} \langle \nabla f , h \rangle_{\mathcal Z} 
= - \int _{\Omega} f (x) \nabla h(x) ~dx , \quad h \in \mathcal Z.
\end{equation}

\noindent {\bf Remark. } 
The definition above applies only to 1st order derivatives. 
If we want to define the second-order derivatives, 
then we need to consider an integral on the boundary. 
This is because our test function space $\mathcal Z$ 
requires that its elements have a value of zero for all 
even-order derivatives $(-\Delta_D)^M f$ for any $M$.

\begin{lem}\label{lem:0806-8} 
{\rm (}Lemma~4.5 in \cite{IMT-2019}{\rm )} 
Let $f \in \mathcal Z'$. Then 
$f = \displaystyle \sum _{j \in \mathbb Z} 
\phi_j(\Lambda _D) f $ in $\mathcal Z'$ 
and $\phi_j(\Lambda _D) f \in L^\infty$ for each $j \in \mathbb Z$. 
\end{lem}

\noindent 
{\bf Definition. } Let $s \in \mathbb R$ and 
$1 \leq p,q \leq \infty$. We define $\dot B^s_{p,q}$ by 
\[
\dot B^s_{p,q} = 
\Big\{ f \in \mathcal Z' \mid
\| f \|_{\dot B^s_{p,q}} = 
\Big\{ \sum _{j \in \mathbb Z} 
    \Big( 2^{sj} \| \phi_j(\Lambda _D) f \|_{L^p} 
    \Big) 
\Big\}^{\frac{1}{q}} 
< \infty 
\Big\} . 
\]

\begin{lem}
Let $s \in \mathbb R$, $1 \leq p \leq \infty$ and $1 \leq q < \infty$. 
Then $\mathcal Z$ is dense in $\dot B^s_{p,q}$. 
\end{lem}

\begin{pf}
Let $q < \infty$ and $f \in \dot B^s_{p,q}$. Then 
$\displaystyle 
f_N = \sum _{|j| \leq N} \phi_j(\Lambda_D)f
$
is an approximation of $f$ 
in the topology of $\dot B^s_{p,q}$, 
and we also see that 
$f_N \in \mathcal Z$ since the domain $\Omega$ is bounded. 
\end{pf}

\begin{lem}{\rm(}\cite{IMT-2019}{\rm)} \label{lem:0806-3}
{\rm (}Dual spaces of Besov spaces{\rm )}
Let $s \in \mathbb R$ and $1 \leq p,q < \infty$. Then 
the dual space of $\dot B^s_{p,q}$ is 
$\dot B^{-s}_{p',q'}$, where $p',q'$ satisfy 
\[
\dfrac{1}{p} + \dfrac{1}{p'} = \dfrac{1}{q} + \dfrac{1}{q'} = 1. 
\]
We also have that 
\[
\| f \|_{\dot B^{-s}_{p',q'}} 
\simeq \sup_{\| g \|_{\dot B^s_{p,q}}=1}  
\big| {}_{\mathcal Z'}\langle f, g \rangle _{\mathcal Z} 
\big| 
= 
\sup_{\| g \|_{\dot B^s_{p,q}}=1}  
\Big| \sum_{j \in \mathbb Z} \int _{\Omega}
  \phi_j(\sqrt{A}) f \cdot 
  \overline{\Phi_j (\sqrt{A})g} ~\mathrm{d}x 
\Big| ,
\]
where $\Phi_j := \phi_{j-1} + \phi_j + \phi_{j+1}$.
\end{lem}

\begin{lem}
If $f,g \in L^2$ and $f = g$ in $\dot B^0_{2,\infty}$, 
then $f = g$ in $L^2$. 
\end{lem}
\begin{pf}
Since $f = g$ in $\dot B^0_{2,\infty}$, 
$\phi_j(\Lambda_D)f = \phi_j(\Lambda_D)g$ in $L^2$. 
We also know the resolution of the identity that 
\[
f = \sum_{j \in \mathbb Z} \phi_j(\Lambda_D)f, 
\quad g = \sum_{j \in \mathbb Z} \phi_j(\Lambda_D)g 
\]
in $L^2$. 
We conclude $f = g$ in $L^2$. 
\end{pf}

We give the following definition for the product of two functions. 

\vskip3mm 

\noindent 
{\bf Definition. } 
Let $f,g \in \mathcal Z'$. Suppose that the following 
functional is continuous. 
\[
\mathcal Z \ni h \mapsto \sum_{k,l \in \mathbb Z} 
{}_{\mathcal Z'} 
\Big\langle \Big( \phi_k(\Lambda_D) f \Big) \Big( \phi_l (\Lambda_D) g \Big) , h 
\Big\rangle _{\mathcal Z} 
\in \mathbb C.
\]
We define $fg$ as an element of $\mathcal Z'$ by 
\[
{}_{\mathcal Z'} \langle fg , h \rangle _{\mathcal Z} 
= 
\sum_{k,l \in \mathbb Z} 
{}_{\mathcal Z'} 
\Big\langle \Big( \phi_k(\Lambda_D) f \Big) \Big( \phi_l (\Lambda_D) g \Big) , h 
\Big\rangle _{\mathcal Z} , 
\quad h \in \mathcal Z.
\]

\subsection{Spectral Restriction and Derivative Estimates}

We introduce several estimates needed in our proof. 
In addition to 
the boundedness of the spectral multiplier from Lemma~\ref{lem:0809-1}, 
we use estimates for all the 2nd derivatives.

\begin{lem}\label{lem:0806-1}
\begin{enumerate}
\item If $f \in L^2$ and $\Lambda _D f \in L^2$, then 
\[
\| \nabla f \|_{L^2} 
= \| \Lambda_Df \|_{L^2}. 
\]
\item If $f \in L^2 $ and $(-\Delta _D) f \in L^2$, then 
\[
\| \nabla^2 f \|_{L^2}  
\leq C \| (-\Delta _D)f \|_{L^2}. 
\]
\item 
If $f, g \in L^2$, then 
\[
\| f \nabla \Lambda_D^{-1} g \|_{L^1} 
\leq C \| f \|_{L^2} \| g \|_{L^2}. 
\]
\end{enumerate}
\end{lem}

\begin{pf}
(1) We use the resolution of the identity that 
$\displaystyle f = \sum _{j \in \mathbb Z} \phi_j(\Lambda _D) f$, 
and its approximation 
$\displaystyle f_N = \sum _{|j| \leq N} \phi_j(\Lambda_D) f$ 
for $N = 1, 2, \dots$. 
Since $(-\Delta_D)f_N \in L^2$, we have 
\[
\| \nabla f_N \|_{L^2}^2 
= \int_{\Omega} f_N (-\Delta _D) f_N \, \mathrm{d}x
= \| \Lambda_D f_N \|_{L^2} ^2.
\]
By taking the limit as $N \to \infty$, 
we obtain $\| \nabla f \|_{L^2} = \| \Lambda_D f \|_{L^2}$. 
\\
(2) The inequality follows from the elliptic estimates on $L^2$ 
(see~\cite{GiTr_2001}). 

\noindent (3) The inequality is obtained by the H\"older inequality and 
 (1). 
\end{pf}

\begin{lem}{\rm(}\cite{FIK-2024, Iw-2018}{\rm)}\label{lem:0806-2} 
Let $1 \leq p \leq \infty$ and $f \in L^p$. 
We have 
\[
\| \nabla e^{t\Delta _D} f \|_{L^p} \leq C t^{-\frac{1}{2}}\| f \|_{L^p}, 
\quad \text{ for all } t> 0,
\]
\[
\| \nabla^2 e^{t\Delta _D} f \|_{L^p} \leq C t^{-1}\| f \|_{L^p}, 
\quad \text{ for all } t> 0,
\]
\[
\| e^{t\Delta _D} \phi_j(\Lambda_D)f \|_{L^p} 
\leq C e^{-c t 2^{2j}}\| \phi_j(\Lambda_D)f \|_{L^p} , 
\quad \text{ for all } t> 0 \text{ and } j \in \mathbb Z. 
\]
\end{lem}

\begin{prop}\label{prop:0809-3}
Let $1 \leq p \leq \infty$. We have 
\[
\sup _{j \in \mathbb Z} 2^{-j}\| \nabla \phi_j (\Lambda _D) \|_{L^p \to L^p}, 
\sup _{j \in \mathbb Z} 2^{-j}\| \phi_j (\Lambda _D) \nabla \|_{L^p \to L^p} 
< \infty .
\]
\[
\sup _{j \in \mathbb Z} 2^{-2j}\| \nabla^2 \phi_j (\Lambda _D) \|_{L^p \to L^p} 
< \infty .
\]
\end{prop}
\begin{pf}
Let $\Phi_j = \phi_{j-1} + \phi_j + \phi_{j+1}$. Then we can write 
\[
\nabla^\alpha \phi_j(\Lambda _D) f
= \nabla^\alpha e^{t\Delta _D} 
  \big(e^{-t\Delta _D} \Phi_j(\Lambda _D)  \big)
  \big( \phi_j(\Lambda_D)f \big), 
\]
where $\alpha$ is a multiindex such that 
$|\alpha |=1,2$. From Lemma~\ref{lem:0806-2} and Lemma~\ref{lem:0809-1}, 
it follows that 
\[
\sup _{j \in \mathbb Z} 2^{-j}\| \nabla \phi_j (\Lambda _D) \|_{L^p \to L^p}, 
\sup _{j \in \mathbb Z} 2^{-2j}\| \nabla^2 \phi_j (\Lambda _D) \|_{L^p \to L^p} 
< \infty .
\]
For $\phi_j(\Lambda_D) \nabla f$, recalling the definition \eqref{0904-8}, 
we reduce the problem to the boundedness 
of $\nabla \phi_j(\Lambda _D)$ on $L^{p'}$ with $1/p + 1/p' = 1$, 
which has already been proved.
\end{pf}

\vskip3mm 

The above and the following two propositions provide elliptic estimates. In this paper, the required inequalities involve derivatives up to the second order, specifically $\nabla^2$.

\begin{prop}{\rm(}\cite{GiTr_2001}{\rm)}\label{prop:0809-4}
Let $1 < p < \infty$. 
If $f$ and $ -\Delta _D f$ belong to $L^p$, then 
\[
\| \nabla^2 f \|_{L^p} 
\leq C \| (-\Delta_D) f \|_{L^p}. 
\]
If $f $  and $ \Lambda _D f $ belong to $L^p$, then 
\[
\| \nabla f \|_{L^p} 
\leq C \| \Lambda_D f \|_{L^p}. 
\]
\end{prop}

\begin{pf}
The inequality for the second-order derivatives is simply an elliptic estimate, 
as proven in \cite{GiTr_2001}. The inequality for the first-order derivatives 
can be established using complex interpolation, as demonstrated in Lemma 4.2 
of \cite{Miy-1980}.
\end{pf}

\section{Proof of Theorem~\ref{Thm:1}}\label{sec:3}

By Lemma~\ref{lem:0806-8}, we can write
$\displaystyle f = \sum_{k \in \mathbb{Z}} f_k$ and 
$\displaystyle g = \sum_{l \in \mathbb{Z}} g_l$, 
where $f_k = \phi_k(\Lambda_D) f$ and 
$g_l = \phi_l(\Lambda_D) g$, with $f_k$ and $g_l$ belonging to $L^\infty$. 
For each $k, l$, it follows from Lemma~\ref{lem:0806-1} and the boundedness of the domain $\Omega$ that 
$(\nabla \Lambda_D^{-1} f_k) g_l$ and 
$f_k \nabla \Lambda_D^{-1} g_l$ 
belong to $L^1$, which is a subspace of $\mathcal{Z}'$.

We prove the following inequality:
\[
\begin{split}
& \Big\{ \sum_{j \in \mathbb{Z}} 
 \Big( 2^{sj}  
    \Big\| 
       \phi_j (\Lambda_D) 
         \sum_{k,l \in \mathbb{Z}} 
      \big( (\nabla \Lambda_D^{-1} f_k) g_l 
           + f_k (\nabla \Lambda_D^{-1} g_l) \big) 
    \Big\|_{L^p} 
  \Big) ^q 
  \Big\}^{\frac{1}{q}} \\
\leq & C \big( \| f \|_{\dot{B}^s_{p_1,q}} \| g \|_{L^{p_2}}  
      + \| f \|_{L^{p_3}} \| g \|_{\dot{B}^s_{p_2,q}} \big).
\end{split}
\]
Once the above inequality is established, 
it follows that 
$(\nabla \Lambda_D^{-1} f) g + f (\nabla \Lambda_D^{-1} g)$ 
belongs to $\mathcal Z'$. 
Therefore, it suffices to demonstrate this inequality. 
To this end, we write
\begin{equation}\label{0820-2}
fg = \sum_{k \in \mathbb{Z}} \sum_{l \leq k} f_k g_l 
 + \sum_{l \in \mathbb{Z}} \sum_{k < l} f_k g_l.
\end{equation}

\noindent 
$\bullet$ \underline{Estimate when $k \geq l$ with $j \geq k$}. 
We mainly explain the estimate for the first term.

We notice that 
$(\nabla \Lambda_D^{-1} f_k) g_l$ 
also satisfies the Dirichlet boundary 
condition and belongs to the domain of the Dirichlet Laplacian. 
The Dirichlet Laplacian can act on the term, and we have 
from Lemma~\ref{lem:0809-1} that for each $j$
\[
\begin{split}
& \Big\| \phi_j (\Lambda_D) 
\sum_{k \leq j} \Big( \big(\nabla \Lambda_D^{-1} f_k \big) 
 \sum_{l \leq k} g_l \Big) \Big\|_{L^p} \\
=& 
 \Big\| 
   \sum_{k \leq j} (-\Delta_D)^{-1} \phi_j (\Lambda_D) (-\Delta_D) 
      \Big( \big(\nabla \Lambda_D^{-1} f_k \big) 
 \sum_{l \leq k} g_l \Big) \Big\|_{L^p} \\
\leq& C 2^{-2j}
\sum_{k \leq j} 
 \Big\| (-\Delta_D) 
 \Big( \big(\nabla \Lambda_D^{-1} f_k \big) 
 \sum_{l \leq k} g_l \Big) \Big\|_{L^p}. 
\end{split}
\]

We apply the Leibniz rule by regarding $-\Delta_D$ 
as the second derivative $-\Delta$ in the weak sense and then estimate 
using the H\"older inequality. For each $k$, we have
\[
\begin{split}
& \| (-\Delta) \nabla \Lambda_D^{-1} f_k \|_{L^{p_1}} 
=  \| \nabla (-\Delta _D)\Lambda _D^{-1} f_k \|_{L^{p_1}}
= \| \nabla \Lambda _D f_k \|_{L^{p_1}}
\leq C 2^{2k} \| f_k \|_{L^{p_1}}, 
\\ 
&\sum_{l \leq k}\| (\nabla^2 \Lambda_{D}^{-1} f_k) \cdot \nabla g_l\|_{L^{p}}
\leq 
 \sum_{l \leq k}\| \nabla ^2 \Lambda_D^{-1} f_k \|_{L^{p_1}}
  2^l \| g_l \|_{L^{p_2}}
  \leq 2^{2k} \| f_k \|_{L^{p_1}} \| g\|_{L^{p_2}}, 
\\
&\sum_{l \leq k} 
\| (\nabla \Lambda_D^{-1} f_k) (-\Delta _D)g_l \|_{L^p}
\leq 
 \sum_{l \leq k} \| \nabla \Lambda _D^{-1} f_k \|_{L^{p_1}}
   2^{2l}\| g_l\|_{L^{p_2}}
\leq C 2^{2k} \| f_k  \|_{L^{p_1}} \| g \|_{L^{p_2}} ,
\end{split}
\]
where we have used 
Proposition~\ref{prop:0809-3}, Lemma~\ref{lem:0809-1} in the first inequality, 
the second derivative estimate with the spectral restriction 
in Proposition~\ref{prop:0809-3} for $\Lambda_D^{-1} f_k$ in the second inequality, 
and the derivative estimates in Proposition~\ref{prop:0809-3} for 
$\Lambda_D^{-1} f_k$ and Lemma~\ref{lem:0809-1} for $g$ in the third inequality. 
Therefore, we obtain
\[
\begin{split}
 2^{-2j}
\sum_{k \leq j} 
 \Big\| (-\Delta _D) 
 \Big( \big(\nabla  \Lambda_D^{-1} f_k \big) 
 \sum_{l \leq k}g_l \Big) \Big\|_{L^p}  
\leq & 
2^{-2j}\sum_{k \leq j} 2^{2k} \| f_k \|_{L^{p_1}} \| g \|_{L^{p_2}}. 
\end{split}
\]
To consider the norm of $\dot B^s_{p,q}$, we multiply the above by $2^{sj}$ 
to take the norm of $\ell^q(\mathbb{Z})$, and write by $k = j + k'$ 
with $k' \leq 0$ that
\begin{equation}\label{0809-5}
\begin{split}
& \Big\{ \sum _{j \in \mathbb Z} 
  \Big( 2^{sj}  
  \Big\| \phi_j (\Lambda _D) 
\sum_{k \leq j} \Big( \big(\nabla  \Lambda_D^{-1} f_k \big) 
 \sum_{l \leq k}g_l \Big) \Big\|_{L^p} 
 \Big) ^q \Big\} ^{\frac{1}{q}} 
\\
\leq 
& C \Big\{ \sum _{j \in \mathbb Z} 
  \Big( 2^{sj}  2^{-2j}
\sum_{k \leq j} 2^{2k} \| f_k \|_{L^{p_1}} 
 \Big)^q \Big\} ^{\frac{1}{q}} 
 \| g \|_{L^{p_2}} 
\\
\leq 
& C \sum_{k' \leq 0} 2^{(-s+2)k'} 
   \Big\{ \sum _{j \in \mathbb Z} 
  \Big( 2^{s(j+k')}  \| f_{j+k'} \|_{L^{p_1}} 
 \Big)^q \Big\} ^{\frac{1}{q}} 
 \| g \|_{L^{p_2}} 
= C \| f \|_{\dot B^s_{p_1,q}}\| g \|_{L^{p_2}}.
\end{split}
\end{equation}
We notice that $ s < 2 $ is needed for the convergence of the sum over $ k' < 0 $.

As for the second term $f_k \nabla \Lambda_D^{-1} \displaystyle \sum_{l \leq k} g_l$ when $j > k$, we apply a similar argument to the above with the estimate of the Riesz transform that 
$\Big\| \nabla \Lambda_D^{-1} \displaystyle \sum_{l \leq k} g_l \Big\|_{L^{p_2}} \leq C \Big\| \displaystyle \sum_{l < k} g_l \Big\|_{L^{p_2}}$, in which we need $1 < p_2 < \infty$ to apply Proposition~\ref{prop:0809-4}. We also know 
$\Big\| \displaystyle \sum_{l \leq k} g_l \Big\|_{L^{p_2}} \leq C \| g \|_{L^{p_2}}$ by the uniformity of the scaled spectral multiplier with the dyadic number $2^k$.

\vskip3mm

\noindent $\bullet$ \underline{Estimate when $k \geq l$ with $k > j$}. 
The case when $k > j$ imposes the restriction $s > -1$, and we derive 
a derivative form of the product. 
We define 
\[
F_k = \Lambda _D^{-1}f_k , \quad 
G_k = \Lambda _D^{-1}g_k
\]
and write 
\begin{equation}\label{0819-2}
\begin{split}
     (\nabla \Lambda _D ^{-1} f_k) g_l
      + f_k \nabla \Lambda_D^{-1} g_l 
=  
&      (\nabla F_k) \Lambda_D G_l
      + (\Lambda _D F_k) \nabla G_l 
\\
=  
&  \nabla \Big( F_k \Lambda_D G_l \Big) 
+ (\Lambda _D F_k) \nabla G_l 
   - F_k \nabla \Lambda _D  G_l . 
\end{split}
\end{equation}

The first term can be handled easily since 
$\nabla \big( F_k \Lambda_D G_l \big) = \nabla \big( (\Lambda _D^{-1}f_k) g_l \big)$, 
and we write, by $k = j+k'$ with $k' > 0$, that
\[
\begin{split}
& \Big\{ \sum _{j \in \mathbb{Z}} 
  \Big( 2^{sj}  
  \Big\| \phi_j (\Lambda _D) 
\sum_{k > j} \nabla \Big( (\Lambda _D^{-1}f_k) \sum_{l \leq k} g_l \Big)
  \Big\|_{L^p} 
 \Big)^q \Big\} ^{\frac{1}{q}} 
\\
\leq 
& C 
 \Big\{ \sum _{j \in \mathbb{Z}} 
  \Big( 2^{(s+1)j}  
\sum_{k > j} 2^{-k} \| f_k \|_{L^{p_1}}
  \Big\| \sum_{l \leq  k}g_l  \Big\|_{L^{p_2}} 
 \Big)^q \Big\} ^{\frac{1}{q}} 
\\
\leq 
& C 
\sum_{k' > 0}2^{-(s+1)k'} 
 \Big\{ \sum _{j \in \mathbb{Z}} 
  \Big( 2^{s(j+k')}  
 \| f_{j+k'} \|_{L^{p_1}}
 \Big)^q \Big\} ^{\frac{1}{q}} 
  \| g  \|_{L^{p_2}} 
\\
= & C \| f \|_{\dot{B}^s_{p_1,q}} \| g \|_{L^{p_2}} .
\end{split}
\]
where we have applied Proposition~\ref{prop:0809-3} to estimate 
$\| \phi_j (\Lambda_D)\nabla \|_{L^p \to L^p}$, 
and we need $s > -1$ for the convergence of the sum over $k' \geq 0$.

For the 2nd and the 3rd terms of the last right-hand side of \eqref{0819-2}, 
we apply the following formula:
\begin{equation}\label{0830-7}
\Lambda _D = c_0 \int_0 ^\infty \mu^{-\frac{3}{2}} 
 \Big( 1 - (1- \mu \Delta _D)^{-1} \Big) {\rm d}\mu,
\end{equation}
where $c_0$ is a positive constant. We then write 
\[
\begin{split}
&   (\Lambda _D F_k) \nabla G_l 
   - F_k \nabla \Lambda _D  G_l
\\
=
& \int_0^\infty \mu^{-\frac{3}{2}} 
 \Big\{ F_k \nabla G_l - F_k \nabla G_l 
  + ( (1-\mu\Delta_D)^{-1} F_k) \nabla G_l
   - F_k \nabla (1-\mu\Delta_D)^{-1}G_l
 \Big\}{\rm d}\mu . 
\end{split}
\]
We find the cancellation for the 1st and the 2nd terms of the integrand above. 
We apply $F_k = (1-\mu\Delta_D)(1-\mu\Delta _D)^{-1} F_k$, 
$G_l = (1-\mu\Delta _D)(1-\mu\Delta _D)^{-1} G_l$ 
to the 3rd and 4th terms in the right-hand side. 
We can write that
\[
\begin{split}
&  (\Lambda _D F_k) \nabla G_l 
   - F_k \nabla \Lambda _D  G_l
\\
=
& \int_0^\infty \mu^{-\frac{3}{2}} 
 \Big\{ 
  \big((1-\mu\Delta_D)^{-1} F_k \big)  
  (-\mu\Delta) \nabla (1-\mu\Delta_D)^{-1} G_l   
\\
& \qquad 
  - \big( (-\mu\Delta _D) (1-\mu\Delta_D)^{-1} F_k \big)  
     \nabla (1-\mu\Delta_D)^{-1} G_l 
 \Big\}{\rm d}\mu . 
\end{split}
\]
Regarding $-\Delta_D$ as the weak derivative, we apply the Leibniz rule that
\begin{equation}\label{0819-1}
\begin{split}
&  (\Lambda _D F_k) \nabla G_l 
   - F_k \nabla \Lambda _D  G_l
\\
=
& \int_0^\infty \mu^{-\frac{3}{2}} 
 \Big\{ 
  (-\mu\Delta)  
  \Big( \big( (1-\mu\Delta_D)^{-1} F_k \big)  
    \nabla  (1-\mu\Delta_D)^{-1} G_l 
  \Big)  
\\
& \qquad 
  + \sum_{m=1}^2 2 \mu \, \partial_{x_m} 
   \Big(
    \big( \partial_{x_m} (1-\mu\Delta_D)^{-1} F_k \big)  
    \nabla (1-\mu\Delta_D)^{-1} G_l 
   \Big) 
 \Big\}{\rm d}\mu . 
\end{split}
\end{equation}
We estimate the $L^p$ norm with the spectral restriction by $\phi_j(\Lambda_D)$, 
and have from 
Proposition~\ref{prop:0809-4} for $G_l$, $1 < p_2 < \infty$,  
that
\[
\begin{split}
& 
\big\| \phi_j (\Lambda_D) \sum_{l \leq k} 
\big((\Lambda _D F_k) \nabla G_l 
   - F_k \nabla \Lambda _D  G_l \big) 
\big\|_{L^p}
\\
\leq 
& C 
\int_0^\infty \mu^{-\frac{1}{2}} 
 \Big\{ 
   2^{2j} \dfrac{1}{1+\mu2^{2k}} 
   \cdot 1 
   + 
   2^{j} \dfrac{2^k}{1+\mu2^{2k}} 
   \cdot 1
 \Big\}
 \| F_k \|_{L^{p_1}} 
 \Big\| \Lambda_D (1-\mu\Delta_D)^{-1} \sum_{l \leq k} G_l \Big\|_{L^{p_2}} 
{\rm d}\mu \\
\leq 
& C 
2^j  \| F_k \|_{L^{p_1}} 
 \Big\| \sum_{l \leq k} \Lambda_D G_l \Big\|_{L^{p_2}} ,
\end{split}
\]
provided that $k > j$. 
We then write by $k = j + k'$ with $k' > 0$ that
\[
\begin{split}
& 
\Big\{ \sum_{j \in \mathbb Z} 
 \Big( 2^{sj} 
    \Big\| \phi_j (\Lambda) \sum_{k > j} \sum_{l \leq k}
           \big( (\nabla ^\perp F_k) \Lambda _D G_l 
           - (\nabla ^{\perp} \Lambda _D F_k )G_l\big) 
    \Big\|_{L^p} \Big)^q
\Big\}^{\frac{1}{q}}
\\
\leq 
& C \Big\{ \sum _{j \in \mathbb Z} 
  \Big(  2^{(s+1)j} \sum _{k > j} \| F_k \| _{L^{p_1}} 
           \Big\| \Lambda_D \sum_{l \leq k} G_l \Big\|_{L^{p_2}} 
  \Big) ^q
  \Big\}^{\frac{1}{q}}
\\
\leq 
& C \sum _{k' > 0} 2^{-(s+1)k'}\Big\{ \sum _{j \in \mathbb Z} 
  \Big(  2^{(s+1)(j+k')}  \| F_{j+k'} \| _{L^{p_1}} 
           \Big\| \Lambda_D \sum_{l \leq j+k'} G_l \Big\|_{L^{p_2}} 
  \Big) ^q
  \Big\}^{\frac{1}{q}}
\\
\leq 
& C \| f \|_{\dot B^s_{p,q}} \| g \|_{L^{p_2}} . 
\end{split}
\]
Here we know the convergence of 
\[
\sum_{k' > 0} 2^{-(s+1)k'}
\]
when $s + 1 > 0$ and 
$\Lambda_D G_l = g_l$.

\vskip3mm 

\noindent 
$\bullet$ \underline{Estimate when $l > k$}. 
This case follows from the same argument as $l \leq k$ 
by replacing $f$ and $g$ with each other. Therefore, we obtain 
$\| f \|_{L^{p_3}} \| g \|_{\dot B^s_{p_4,q}}$ 
in the inequality on the right-hand side.

\section{Proof of Theorem~\ref{Thm:2}}\label{sec:4}

We apply the argument from \cite{IwUe-2024}, adapting the spaces 
to those defined on a bounded domain. We begin with the following proposition.

\begin{prop}\label{prop:0806-5}  
Let $1 < p < 2$. Suppose that 
$\theta $ is a mild solution of \eqref{eq:SQG}. 
Then we have 
\[
\| \theta (t) - e^{t\Delta _D} \theta_0 ||_{\dot B^{-1+\frac{2}{p}}_{p,\infty}}
\leq C \| \theta \|_{L^\infty (0,T ; L^2)}^2, 
\quad \text{ for } t \in [0,T]. 
\] 
\end{prop}
\begin{pf}
We have from the definition of the mild solution that for every $h \in \mathcal Z$ 
\[
\Big| 
{}_{\mathcal Z'} \langle \theta - e^{t\Delta _D} \theta_0, h \rangle_{\mathcal Z}
\Big| 
\leq  \int_0^t \Big| \int_{\Omega} u \theta \cdot \nabla e^{(t-\tau)\Delta _D} h
~{\rm d}x \Big| {\rm d}\tau . 
\]
We see that $u \theta \in L^1$ by the H\"older inequality and 
Lemma~\ref{lem:0806-1}. 
It follows from Lemma~\ref{lem:0806-2} and the resolution of the 
identity in Lemma~\ref{lem:0806-8} that 
\[
\begin{split}
& 
\int_0^t \Big| \int_{\Omega} u \theta \cdot \nabla e^{(t-\tau)\Delta _D}h 
~{\rm d}x \Big| {\rm d}\tau
\\
\leq
&  \int_0^t \| u \theta \| _{L^1} 
\Big\| \nabla e^{(t-\tau) \Delta _D} 
      \sum_{j \in \mathbb Z} \phi_j(\Lambda_D) h \Big\|_{L^\infty}
 {\rm d}\tau 
\\
\leq
&  C \| \theta \|_{L^\infty (0,T; L^2)}^2 
\sum_{j \in \mathbb Z} 
\int_0^t (t-\tau)^{-\frac{1}{2}} e^{-c(t-\tau) 2^{2j}}
 {\rm d}\tau 
\Big\| \phi_j(\Lambda_D) h \Big\|_{L^\infty}.
\end{split}
\]
Since $\displaystyle \int_0^t (t-\tau)^{-\frac{1}{2}} e^{-c(t-\tau) 2^{2j}}
 {\rm d}\tau  
 \leq C 2^{-j}$ and 
 $\dot B^{-1+\frac{2}{p'}}_{p',1} \hookrightarrow \dot B^{-1}_{\infty,1}$, 
we have 
\[
\int_0^t \Big| \int_{\Omega} u \theta \cdot \nabla e^{(t-\tau)\Delta _D}h 
~{\rm d}x \Big| {\rm d}\tau
\leq  C
 \| \theta \|_{L^\infty (0,T; L^2)}^2 
 \| h \|_{\dot B^{-1}_{\infty ,1}}
\leq  C
 \| \theta \|_{L^\infty (0,T; L^2)}^2 
 \| h \|_{\dot B^{-1+\frac{2}{p'}}_{p' ,1}},    
\]
where $1/p + 1/p' = 1$. 
We know from Lemma~\ref{lem:0806-3} that 
$\dot B^{-1+\frac{2}{p}}_{p,\infty }$ 
is the topological dual of $ \dot B^{-1+\frac{2}{p'}}_{p',1}
= \dot B^{1- \frac{2}{p}}_{p',1}$, 
and that $\mathcal Z$ is dense in $\dot B^{-1+\frac{2}{p'}}_{p' ,1}$. 
Taking the supremum with respect to $h$ with 
$\| h \|_{\dot B^{-1+\frac{2}{p'}}_{p',1}} = 1$, we obtain 
\[
\| \theta - e^{t \Delta_D} \theta_0 \|_{\dot B^{-1+\frac{2}{p}}_{p,\infty}} 
\leq C \| \theta \|_{L^\infty (0,T ; L^2)} ^2 ,
\]
which completes the proof. 
\end{pf}

\vskip3mm

\noindent 
{\bf Proof of Theorem~\ref{Thm:2}. }  
Let $\theta_0 = \theta(0) = \widetilde \theta(0)$, 
$\psi := \theta - \widetilde \theta$ and $1 < p < 2$. 
We show that $\psi (t) = 0$ in $\dot B^{-1+\frac{2}{p}}_{p,\infty}$ 
for $t \in [0,\delta]$ for some small $\delta > 0$, 
which implies $\psi = 0$ in $L^2$ for $t \in [0,\delta]$.  

We will choose $\delta > 0$ such that the following holds: 
\begin{equation}\label{0826-1}
\displaystyle 
\sup_{t \in (0,\delta]} t^{\frac{1}{2}-\frac{1}{2p}} 
\| e^{t\Delta_D} \theta_0 \|_{L^{2p}}, 
\sup_{t \in [0,\delta]}
\| \theta (t) - e^{t\Delta_D}\theta_0 \|_{L^2}, 
\sup_{t \in [0,\delta]}
\| \widetilde \theta (t) - e^{t\Delta_D} \widetilde \theta_0 \|_{L^2} 
\ll 1. 
\end{equation}
The 2nd and 3rd smallness hold due to 
the assumption of the continuity of the solutions. 
As for the 1st smallness condition, we can analogously prove 
to the entire space case that
\begin{equation}\label{0904-1}
\lim_{T\to 0}\sup_{t\in(0,T]}
t^{\frac{1}{2}-\frac{1}{2p}}  
 \| e^{t\Delta _D}u_{0} \|_{L^{2p}} = 0 . 
\end{equation}
The proof of this is achieved by using the boundedness:
\[
t^{\frac{1}{2}-\frac{1}{2p}}  
 \| e^{t\Delta _D}u_{0} \|_{L^{2p}}
  \leq C \| u_0 \|_{\dot B^{-1+\frac{1}{p}}_{2p,\infty}}
  \leq C \| u_0 \|_{L^2}, 
\]
and the following approximation of the initial data:
\[
u_{0,N} = \sum_{|j| \leq N} \phi_j(\sqrt{A})u_0.
\]
Since $u_{0,N}$ satisfies 
$(-\Delta _D)^M u_{0,N}$ for all $M \in \mathbb N$, 
we have 
\[
\lim_{T\to 0}\sup_{t\in(0,T]}
t^{\frac{1}{2}-\frac{1}{2p}}  
 \| e^{t\Delta _D}u_{0,N} \|_{L^{2p}} = 0 , 
 \quad \text{for each } N.  
\]
We also see that $u_{0,N}$ converges to $u_0$ in 
$L^2$. 
Thus, the proof of \eqref{0904-1} is complete.

In what follows we investigate the nonlinear estimates 
and prove the uniqueness, dividing it into several steps.

\vskip3mm 

\noindent 
{\bf Step 1. } 
We write 
\[
\psi = -\dfrac{1}{2}\int_0^t e^{(t-\tau)\Delta _D}
\nabla \cdot 
\Big\{ (\nabla^\perp \Lambda_D^{-1} \psi) \theta 
          + (\nabla ^\perp \Lambda _D^{-1} \theta ) \psi 
       + (\nabla^\perp \Lambda_D^{-1} \psi) \widetilde \theta 
          + (\nabla ^\perp \Lambda _D^{-1} \widetilde \theta \, ) \psi 
\Big\}
{\rm d}\tau .
\]
We only consider the first two terms, as the 3rd and 4th terms 
are handled in an identical manner. 
We write 
\[
N(\theta) := \theta - e^{\tau \Delta_D}\theta _0,
\]
and 
\[
\begin{split}
I:= &
\Big\| 
\int_0^t e^{(t-\tau)\Delta _D}
\nabla \cdot 
\Big\{ (\nabla^\perp \Lambda_D^{-1} \psi) \theta 
          + (\nabla ^\perp \Lambda _D^{-1} \theta ) \psi 
\Big\}
{\rm d}\tau 
\Big\|_{\dot B^{-1+\frac{2}{p}}_{p,\infty}} 
\\
\leq 
& 
\Big\| 
\int_0^t e^{(t-\tau)\Delta _D}
\nabla \cdot 
\Big\{ (\nabla^\perp \Lambda_D^{-1} \psi) e^{\tau \Delta_D}\theta_0
          + \big(\nabla ^\perp \Lambda _D^{-1} e^{\tau \Delta_D}\theta_0 \big) \psi 
\Big\} 
{\rm d}\tau 
\Big\|_{\dot B^{-1+\frac{2}{p}}_{p,\infty}} 
\\
& 
+ 
\Big\| 
\int_0^t e^{(t-\tau)\Delta _D}
\nabla \cdot 
\Big\{ (\nabla^\perp \Lambda_D^{-1} \psi) N(\theta )
          + \big(\nabla ^\perp \Lambda _D^{-1} N(\theta) \big) \psi  
\Big\}
{\rm d}\tau 
\Big\|_{\dot B^{-1+\frac{2}{p}}_{p,\infty}}  
\\
=: & I_{\rm linear} + I_{\rm nonlinear} .
\end{split}
\]

\noindent 
{\bf Step 2. } 
We show the following. 
\begin{equation}\label{0820-5}
 I_{\rm linear}
 \leq 
C 
\Big( \sup_{\tau \in [0,\delta]} \|\psi\|_{\dot B^{-1+\frac{2}{p}}_{p,\infty}} 
\Big) 
\Big( \sup_{\tau \in (0,\delta]} \tau ^{\frac{1}{2} - \frac{1}{2p}} 
     \| e^{\tau \Delta _D}\theta_ 0\|_{L^{2p}}  
\Big)  . 
\end{equation}
By Lemma~\ref{lem:0806-3}, we consider the coupling 
with $h \in \dot B^{1-\frac{2}{p}}_{p,1}$ 
\[
\begin{split}
& \int_{\Omega} 
\Big[ \int_0^t e^{(t-\tau)\Delta _D}
\nabla \cdot 
\Big\{ (\nabla^\perp \Lambda_D^{-1} \psi) e^{\tau \Delta_D}\theta_0
          + \big(\nabla ^\perp \Lambda _D^{-1} e^{\tau \Delta_D}\theta_0 \big) \psi 
\Big\} 
{\rm d}\tau 
\Big] 
h ~{\rm d}x, 
\end{split}
\]
instead of the original norm in $\dot B^{-1+\frac{2}{p}}_{p,\infty}$.   
As in the decomposition \eqref{0820-2}, 
we divide into the two cases: $k \geq l$ and $k < l$, and write
\[
\begin{split}
I_{\rm linear } ^1
= 
&
\sum_{k\geq l} 
\int_{\Omega} 
\Big[ \int_0^t e^{(t-\tau)\Delta _D}
\nabla \cdot 
\Big\{ (\nabla^\perp \Lambda_D^{-1} \psi_k) (e^{\tau \Delta_D}\theta_0)_l
          + \big(\nabla ^\perp \Lambda _D^{-1} (e^{\tau \Delta_D}\theta_0)_l \big) \psi_k 
\Big\} 
{\rm d}\tau 
\Big] 
h ~{\rm d}x, 
\\
I_{\rm linear } ^2
= 
&
\sum_{k< l} 
\int_{\Omega} 
\Big[ \int_0^t e^{(t-\tau)\Delta _D}
\nabla \cdot 
\Big\{ (\nabla^\perp \Lambda_D^{-1} \psi_k) (e^{\tau \Delta_D}\theta_0)_l
          + \big(\nabla ^\perp \Lambda _D^{-1} (e^{\tau \Delta_D}\theta_0)_l \big) \psi_k 
\Big\} 
{\rm d}\tau 
\Big] 
h ~{\rm d}x.
\end{split}
\]
We can also consider the dual operator of $e^{(t-\tau)\Delta_D} \nabla \cdot$ 
and the decomposition $h = \displaystyle \sum_{j \in \mathbb{Z}} h_j$, 
where $h_j = \phi_j(\Lambda_D) h$.

\vskip2mm  

\noindent 
Step 2.1. \,  
If we consider the sum with the restriction $j \geq k$ 
in $I_{\rm linear}^1$, then  
\begin{equation}\label{0820-21}
\begin{split}
& 
I_{{\rm linear}, \, j \geq k}^1 
\\
=
& 
\Big| 
\sum_{j \geq k\geq l} 
\int_{\Omega} 
\Big[ \int_0^t e^{(t-\tau)\Delta _D}
\nabla \cdot 
\Big\{ (\nabla^\perp \Lambda_D^{-1} \psi_k) (e^{\tau \Delta_D}\theta_0)_l
          + \big(\nabla ^\perp \Lambda _D^{-1} (e^{\tau \Delta_D}\theta_0)_l \big) \psi_k 
\Big\} 
{\rm d}\tau 
\Big] 
h_j ~{\rm d}x
\Big| 
\\
= & 
\Big| 
\sum_{j  \geq k} 
\int_0^t
\int_{\Omega} 
\Big\{ (\nabla^\perp \Lambda_D^{-1} \psi_k) 
  \sum_{l \leq k}(e^{\tau \Delta_D}\theta_0)_l
          + \Big(\nabla ^\perp \Lambda _D^{-1} \sum_{l \leq k}(e^{\tau \Delta_D}\theta_0)_l \Big) \psi_k 
\Big\} 
\cdot \nabla e^{(t-\tau)\Delta _D} h_j ~{\rm d}x
~{\rm d}\tau 
\Big| 
\\
\leq  & 
\sum_{j \in \mathbb Z } \sum _{k \leq j} 
\int_0^t
\| \psi_k \|_{L^{2p}}  \|e^{\tau \Delta_D}\theta_0\|_{L^{2p}}
\| \nabla e^{(t-\tau)\Delta _D} h_j \|_{L^{p'}}~{\rm d}x
~{\rm d}\tau ,
\end{split}
\end{equation}
where $1/p + 1/p' = 1$. Lemma~\ref{lem:0809-1} yields that 
\[
\sum_{k \leq j} 
\| \psi_k \|_{L^{2p}} 
\leq C \sum _{k \leq j} 2^{2(\frac{1}{2}-\frac{1}{2p})k} \| \psi_k \|_{L^2}
\leq C 2^{(1-\frac{1}{p})j}\| \psi \|_{\dot B^0_{2,\infty}},
\]
which implies that 
\begin{equation}\label{0822-2}
\begin{split}
& 
I_{{\rm linear}, \, j \geq k}^1 
\\
\leq 
& C \sum_{j \in \mathbb Z}  \int _0^t 2^{(1-\frac{1}{p})j}
\cdot \tau ^{-\frac{1}{2}+\frac{1}{2p}}
  \cdot 2^j \cdot 
  (t-\tau) ^{-\frac{1}{2}(1+\frac{1}{p})} \cdot 2^{-(1+\frac{1}{p})j} 
  {\rm d}\tau 
\\
&\qquad 
  \cdot 
  \Big( \sup_{\tau (0,\delta]}\| \psi \|_{\dot B^{0}_{2,\infty}} 
  \Big)
  \Big( \sup_{\tau\in(0,\delta]} \tau^{\frac{1}{2}-\frac{1}{2p}} 
     \| e^{\tau \Delta _D} \theta _0 \|_{L^{2p}} 
  \Big) 
  \| h_j \|_{L^{p'}}
\\
= 
& C 
\Big( \int_0^1 \tau ^{-\frac{1}{2}+\frac{1}{2p}}
   (1-\tau)^{-\frac{1}{2}-\frac{1}{2p}} {\rm d}\tau 
\Big) 
  \Big( \sup_{\tau (0,\delta]}\| \psi \|_{\dot B^{0}_{2,\infty}} 
  \Big)
  \Big( \sup_{\tau\in(0,\delta]} \tau^{\frac{1}{2}-\frac{1}{2p}} 
     \| e^{\tau \Delta _D} \theta _0 \|_{L^{2p}} 
  \Big) 
  \| h \|_{\dot B^{1-\frac{2}{p}}_{p',\infty}} .
\end{split}
\end{equation}
This proves a part of \eqref{0820-5} by taking 
the supremum over $h$ such that 
$\| h \|_{\dot B^{1-\frac{2}{p}}_{p,\infty}}=1$.

For the sum with the restriction $j < k$, we need an analogous argument 
to \eqref{0819-1} with $\nabla ^\perp$, instead of $\nabla $. 
We need to consider 
\[
\begin{split}
\nabla \cdot 
\Big\{ (\nabla^\perp \Lambda_D^{-1} \psi_k) (e^{\tau \Delta_D}\theta_0)_l
          + \big(\nabla ^\perp \Lambda _D^{-1} (e^{\tau \Delta_D}\theta_0)_l \big) \psi_k 
\Big\} .
%& 
%\int_{\Omega} e^{(t-\tau)\Delta _D}
%\nabla \cdot 
%\Big\{ (\nabla^\perp \Lambda_D^{-1} \psi_k) (e^{\tau \Delta_D}\theta_0)_l
%          + \big(\nabla ^\perp \Lambda _D^{-1} (e^{\tau \Delta_D}\theta_0)_l \big) \psi_k 
%\Big\} 
%h_j~
%{\rm d}x
%\\
%=& 
%\int_{\Omega} 
%\Big\{ (\nabla^\perp \Lambda_D^{-1} \psi_k) (e^{\tau \Delta_D}\theta_0)_l
%          + \big(\nabla ^\perp \Lambda _D^{-1} (e^{\tau \Delta_D}\theta_0)_l \big) \psi_k 
%\Big\} 
%\cdot \nabla e^{(t-\tau)\Delta _D}h_j~
%{\rm d}x ,
\end{split}
\]
%since $ e^{(t-\tau)\Delta _D}h_j$ satisfies the Dirichlet boundary condition 
%and integration by parts can be applied here. 
The corresponding first term on the last right-hand side of 
\eqref{0819-2} is 
$\nabla ^\perp \Big( (\Lambda _D^{-1}\psi_k) (e^{\tau \Delta _D}\theta_0)_l 
\Big)$ and 
\[
\nabla \cdot \nabla ^\perp 
\Big( (\Lambda _D^{-1}\psi_k) (e^{\tau \Delta _D}\theta_0)_l 
\Big)
= 0.
\]
For the corresponding 2nd and 3rd terms to \eqref{0819-2}, we write 
\begin{equation}\label{0822-5}
F_k = \Lambda _D^{-1} \psi_k, 
\qquad 
G_l = \Lambda _D^{-1} (e^{\tau \Delta _D}\theta_0)_l. 
\end{equation}
We rewrite the equality \eqref{0819-1} 
involving $\nabla^\perp$ instead of $\nabla$, 
by the Leibniz rule,  
and then need to consider the following. 
\begin{equation}\notag 
\begin{split}
&  (\Lambda _D F_k) \nabla^\perp G_l 
   - F_k \nabla^\perp \Lambda _D  G_l\\
=
& \int_0^\infty \mu^{-\frac{3}{2}} 
\sum_{m=1}^2 
 \Big\{ 
  \mu \Delta 
  \Big( \big( (1-\mu\Delta_D)^{-1} F_k \big)  
   \nabla^\perp (1-\mu\Delta_D)^{-1} G_l 
  \Big)  
\\
& \qquad 
   - 2\mu \, \partial_{x_m} 
   \Big(
    \big( (1-\mu\Delta_D)^{-1} F_k \big)  
     \partial _{x_m}\nabla^\perp (1-\mu\Delta_D)^{-1} G_l 
   \Big) 
 \Big\}{\rm d}\mu . 
\end{split}
\end{equation}
Its divergence is written using the Leibniz rule and 
the fact that $\nabla \cdot \nabla ^\perp = 0$ as follows. 
\begin{equation}\notag 
\begin{split}
&  \nabla \cdot 
\Big( (\Lambda _D F_k) \nabla^\perp G_l 
   - F_k \nabla^\perp \Lambda _D  G_l
\Big) 
\\
=
& \int_0^\infty \mu^{-\frac{1}{2}} 
\sum_{m=1}^2 
 \Big\{ 
  \Delta_D 
  \Big( \big( \nabla (1-\mu\Delta_D)^{-1} F_k \big)\cdot   
   \nabla^\perp (1-\mu\Delta_D)^{-1} G_l 
  \Big)  
\\
& \qquad 
   - 2 \nabla \cdot \partial_{x_m} 
   \Big(
    \big( (1-\mu\Delta_D)^{-1} F_k \big)  
     \partial _{x_m}\nabla^\perp (1-\mu\Delta_D)^{-1} G_l 
   \Big) 
 \Big\}{\rm d}\mu , 
\end{split}
\end{equation}
where the derivative $\mu \Delta $ in the first integrand 
is the Dirichlet Laplacian $\Delta _D$. 
This is because 
\[
\big( \nabla (1-\mu\Delta_D)^{-1} F_k \big)\cdot   
   \nabla^\perp (1-\mu\Delta_D)^{-1} G_l 
\]
satisfies  the Dirichlet boundary condition 
   (see Lemma~3.3 in \cite{CoNg-2018-2}). 
For the second term, we note that 
\[
\big( (1-\mu\Delta_D)^{-1} F_k \big)  
    \partial _{x_m} \nabla ^\perp (1-\mu\Delta_D)^{-1} G_l 
\]
satisfy the Dirichlet boundary condition. 
Since $h_j$ also satisfies the Dirichlet boundary condition, 
we apply integration by parts twice, which implies that 
\begin{equation}\label{0823-2}
\begin{split}
& 
\int_{\Omega} e^{(t-\tau)\Delta _D}
\Big\{ \nabla \cdot 
\Big( (\nabla^\perp \Lambda_D^{-1} \psi_k) (e^{\tau \Delta_D}\theta_0)_l
          + \big(\nabla ^\perp \Lambda _D^{-1} (e^{\tau \Delta_D}\theta_0)_l \big) \psi_k 
\Big) 
\Big\} 
h_j~
{\rm d}x
\\
=& 
\sum_{m=1}^2 \int_0^\infty \mu^{-\frac{1}{2}} \int_{\Omega} 
 \Big\{ 
  \Big( \big( \nabla (1-\mu\Delta_D)^{-1} F_k \big)\cdot   
   \nabla^\perp (1-\mu\Delta_D)^{-1} G_l 
  \Big)  
  \Delta _D  e^{(t-\tau)\Delta _D}h_j
\\
& \qquad 
   - 
   2\Big(
    \big( (1-\mu\Delta_D)^{-1} F_k \big)  
     \partial _{x_m}\nabla^\perp (1-\mu\Delta_D)^{-1} G_l 
   \Big) 
   \cdot \nabla  \partial_{x_m} e^{(t-\tau)\Delta _D}h_j
 \Big\}  
  ~{\rm d}x {\rm d}\mu . 
\end{split}
\end{equation}
We then obtain by the equality above that 
\[
\begin{split}
& I_{{\rm linear}, \, j < k}^1 
\\
=& 
\Big| \sum_{j \in \mathbb Z} \sum_{k > j} \sum _{l \leq k}
\int_0^t  \int_{\Omega} 
\Big[ e^{(t-\tau)\Delta _D}
\nabla \cdot 
\Big\{ (\nabla^\perp \Lambda_D^{-1} \psi_k) (e^{\tau \Delta_D}\theta_0)_l
          + \big(\nabla ^\perp \Lambda _D^{-1} (e^{\tau \Delta_D}\theta_0)_l \big) \psi_k 
\Big\}
\Big] 
\\
& \qquad  
h_j~
{\rm d}x {\rm d}\tau  
\Big| 
\\
\leq & C  
\sum_{j \in \mathbb Z} \sum_{k > j} 
\int_0^t 
\int_0^\infty \mu^{-\frac{1}{2}} 
\cdot  \dfrac{1}{1+\mu 2^{2k}} \| F_k \|_{L^{2p}}
  \sum _{l \leq k} 2^{ 2l}
 \| G_l \|_{L^{2p}}
\\
  & \qquad 
\cdot 2^{2j} (t-\tau)^{-\frac{1}{2}(1+\frac{1}{p})} 
\cdot 2^{-(1+\frac{1}{p})j} 
\| h_j \|_{L^{p'}}
~{\rm d}\mu {\rm d}\tau ,
\end{split}
\]
where $1/p + 1/p' = 1$. 
We change the variables 
$\mu $ and $ \tau$ with $2^{-2k}\mu $ and $ t\tau$, 
respectively, and obtain that 
\begin{equation}\label{0822-1}
\begin{split}
& I_{{\rm linear}, \, j < k}^1 
\\
\leq & C 
\Big( \int_0^\infty \dfrac{\mu^{-\frac{1}{2}}}{1+\mu} {\rm d}\mu \Big)
\Big( \int_0^1 \tau ^{-\frac{1}{2}+\frac{1}{2p}} 
  (1-\tau)^{-\frac{1}{2}-\frac{1}{2p}} {\rm d}\tau \Big)
\\
& 
\cdot 
  \sup_{ \tau \in (0,t)} \tau ^{\frac{1}{2}-\frac{1}{2p}} 
\sum_{j \in \mathbb Z} \sum_{k > j} 
 2^{-k} \| F_k \|_{L^{2p}}
  \sum _{l \leq k} 2^{ 2l}
 \| G_l \|_{L^{2p}}
\cdot 2^{2j} 
\cdot 2^{-(1+\frac{1}{p})j} 
\| h_j \|_{L^{p'}}
\\
\leq & C 
\sup_{ \tau \in (0,t)} \tau ^{\frac{1}{2}-\frac{1}{2p}} 
\Big( 
\sup_{j \in \mathbb Z} \sum_{k > j}  
\| \Lambda _D^{-1} \psi_k \|_{L^{2p}}
 \| e^{\tau \Delta _D} \theta_0 \|_{L^{2p}}
\cdot 2^{\frac{1}{p}j} 
\Big) \Big( \sum_{j \in \mathbb Z}2^{(1-\frac{2}{p})j}\| h_j \|_{L^{p'}}\Big).
\end{split}
\end{equation}
We can write $k = j+ k'$ with $k' > 0$ 
and use the embedding $\dot B^{-1+\frac{1}{p}}_{2p,\infty} 
\hookrightarrow \dot B^{-1+\frac{2}{p}}_{p,\infty}$ 
to obtain that 
\[
\begin{split}
\sup_{j \in \mathbb Z} \sum_{k > j}  
\| \Lambda _D^{-1} \psi_k \|_{L^{2p}}
\cdot 2^{\frac{1}{p}j} 
\leq
& C \sum_{k' > 0} 2^{-\frac{1}{p}k'} 
\sup_{j \in \mathbb Z} 
2^{(-1 +\frac{1}{p})(j+k')} \| \psi_{j+k'} \|_{L^{2p}}
\leq 
C \| \psi \|_{\dot B^{-1+\frac{2}{p}}_{p,\infty}},
\end{split}
\]
which implies that 
\[
I_{{\rm linear}, \, j < k}^1 
\leq 
C \Big(\sup _{\tau \in[0,\delta]}\| \psi \|_{\dot B^{-1+\frac{2}{p}}_{p,\infty}} \Big)
\Big( \sup_{\tau \in (0,\delta)} \tau ^{\frac{1}{2}-\frac{1}{2p}} 
   \| e^{\tau \Delta _D} \theta_0 \|_{L^{2p}}
\Big) 
\| h \|_{\dot B^{1- \frac{2}{p}}_{p',1}}. 
\]
The supremum over $h$ such that $\| h \|_{\dot B^{1- \frac{2}{p}}_{p',1}} = 1$ 
proves a part of \eqref{0820-5}.

By the two argument above for 
$I_{{\rm linear}, \, j \geq k}^1 $ and $ I_{{\rm linear}, \, j < k}^1 $, 
we obtain the following inequality for $I^1_{\rm linear}$. 
\[
\sup _{\| h \|_{\dot B^{1-\frac{2}{p}}_{p,\infty}}=1 } I^1_{\rm inear}
\leq C \Big(\sup _{\tau \in[0,\delta]}\| \psi \|_{\dot B^{-1+\frac{2}{p}}_{p,\infty}} \Big)
\Big( \sup_{\tau \in (0,\delta)} \tau ^{\frac{1}{2}-\frac{1}{2p}} 
   \| e^{\tau \Delta _D} \theta_0 \|_{L^{2p}}
\Big).  
\]

\vskip2mm  

\noindent 
Step 2.2. \,  
We now consider $I^2_{\rm linear}$ 
with the coupling involving $h$ 
such that $\|h \|_{\dot B^{1-\frac{2}{p}}_{p,\infty}} = 1$. 
Since $l > k$, we need to study two cases: $j \geq l$ and $j < l$. 

In the case when $j \geq l$, 
\begin{equation}\label{0820-21}
\begin{split}
& 
I_{{\rm linear}, \, j \geq l}^1 
\\
=
& 
\Big| 
\sum_{j \geq l> k} 
\int_{\Omega} 
\Big[ \int_0^t e^{(t-\tau)\Delta _D}
\nabla \cdot 
\Big\{ (\nabla^\perp \Lambda_D^{-1} \psi_k) (e^{\tau \Delta_D}\theta_0)_l
          + \big(\nabla ^\perp \Lambda _D^{-1} (e^{\tau \Delta_D}\theta_0)_l \big) \psi_k 
\Big\} 
{\rm d}\tau 
\Big] 
h_j ~{\rm d}x
\Big| 
\\
\leq  & 
\sum_{j \in \mathbb Z } \sum _{l \leq j} 
\sum_{k < l}
\int_0^t
\| \psi_k \|_{L^{2p}}  
\| (e^{\tau \Delta_D}\theta_0)_l \|_{L^{2p}}
\| \nabla e^{(t-\tau)\Delta _D} h_j \|_{L^{p'}}~{\rm d}x
~{\rm d}\tau ,
\end{split}
\end{equation}
where $1/p + 1/p' = 1$. Lemma~\ref{lem:0809-1} implies that 
\[
\begin{split}
\sum_{l \leq j}\sum_{k <l} 
\| \psi_k \|_{L^{2p}} 
\| (e^{\tau \Delta _D} \theta_0)_l \|_{L^{2p}}
\leq 
& C \| \psi \|_{\dot B^0_{2,\infty}} 
\sum_{l \leq j} 2^{(1-\frac{1}{p})l}
 \| (e^{\tau \Delta _D} \theta_0)_l \|_{L^{2p}}
\\
\leq 
& 2^{(1-\frac{1}{p})j} 
 \| \psi \|_{\dot B^0_{2,\infty}}
 \| e^{\tau \Delta _D} \theta_0 \|_{L^{2p}},
\end{split}
\]
which implies, similarly to \eqref{0822-2}, that 
\begin{equation}\notag 
\begin{split}
I_{{\rm linear}, \, j \geq l}^2 
\leq 
& C 
  \Big( \sup_{\tau \in (0,\delta]}\| \psi \|_{\dot B^{0}_{2,\infty}} 
  \Big)
  \Big( \sup_{\tau\in(0,\delta]} \tau^{\frac{1}{2}-\frac{1}{2p}} 
     \| e^{\tau \Delta _D} \theta _0 \|_{L^{2p}} 
  \Big) 
  \| h \|_{\dot B^{1-\frac{2}{p}}_{p',\infty}} .
\end{split}
\end{equation}

In the case when $j < l$, we write 
\begin{equation}\notag 
\widetilde F_l = \Lambda _D^{-1} (e^{\tau \Delta _D}\theta_0)_l, 
\qquad 
\widetilde G_k = \Lambda _D^{-1} \psi_k, 
\end{equation}
instead of \eqref{0822-5}. 
We apply the same argument as in the first inequality of \eqref{0822-1} 
and estimate as follows. 
\begin{equation}\notag 
\begin{split}
& I_{{\rm linear}, \, j < l}^2 
\\
\leq & C 
  \sup_{ \tau \in (0,t)} \tau ^{\frac{1}{2}-\frac{1}{2p}} 
\sum_{j \in \mathbb Z} \sum_{l > j} 
\Big\{ 
 2^{-k} \| \widetilde F_l \|_{L^{2p}}
  \sum _{k \leq l} 2^{ 2k}
 \| \widetilde G_k \|_{L^{2p}}
+  \| \widetilde F_l \|_{L^{2p}} 
  \Big\| \Lambda_D \sum _{k \leq l} \widetilde G_k \Big\|_{L^{2p}} 
\Big\} 
\\
& 
\qquad \qquad \qquad \qquad \qquad 
\cdot 2^{2j} 
\cdot 2^{-(1+\frac{1}{p})j} 
\| h_j \|_{L^{p'}}
\\
\leq & C 
\sup_{ \tau \in (0,t)} \tau ^{\frac{1}{2}-\frac{1}{2p}} 
\Big( 
\sup_{j \in \mathbb Z} \sum_{l > j}
 \| (e^{\tau \Delta _D} \theta_0)_l \|_{L^{2p}}
 \sum_{k <l}  \| \Lambda _D^{-1} \psi_k \|_{L^{2p}}
\cdot 2^{\frac{1}{p}j} 
\Big) \Big( \sum_{j \in \mathbb Z}2^{(1-\frac{2}{p})j}\| h_j \|_{L^{p'}}\Big)
\\
\leq 
& C 
  \Big( \sup_{\tau\in(0,\delta]} \tau^{\frac{1}{2}-\frac{1}{2p}} 
     \| e^{\tau \Delta _D} \theta _0 \|_{L^{2p}} 
  \Big) 
  \Big( \sup_{\tau \in (0,\delta]}\| \psi \|_{\dot B^{0}_{2,\infty}} 
  \Big)
  \| h \|_{\dot B^{1-\frac{2}{p}}_{p',\infty}} .
\end{split}
\end{equation}

We obtain \eqref{0820-5} by Step 2.1 and Step 2.2.

\vskip3mm 

\noindent 
{\bf Step 3. } 
We show the following. 
\begin{equation}\label{0823-1}
 I_{\rm nonlinear}
 \leq 
C 
\Big( \sup_{\tau \in [0,\delta]} \|\psi\|_{\dot B^{-1+\frac{2}{p}}_{p,\infty}} 
\Big) 
\Big( \sup_{\tau \in (0,\delta]} \tau ^{\frac{1}{2} - \frac{1}{2p}} 
     \| N(\theta )\|_{L^{2}}  
\Big)  . 
\end{equation}
Let $h \in \dot B^{1-\frac{2}{p}}_{p,1}$. By Lemma~\ref{lem:0806-3}, we 
again consider the following coupling: 
\[
\begin{split}
& \int_{\Omega} 
\Big[ \int_0^t e^{(t-\tau)\Delta _D}
\nabla \cdot 
\Big\{ (\nabla^\perp \Lambda_D^{-1} \psi) N(\theta)
          + \big(\nabla ^\perp \Lambda _D^{-1} N(\theta) \big) \psi 
\Big\} 
{\rm d}\tau 
\Big] 
h ~{\rm d}x 
\\
= & 
- \sum_{j\in \mathbb Z} 
 \Big( \sum _{k \geq l} + \sum _{k <l} \Big) 
\int_0^t 
\int_{\Omega} 
\Big[  
\Big\{ (\nabla^\perp \Lambda_D^{-1} \psi_k) N(\theta)_l
          + \big(\nabla ^\perp \Lambda _D^{-1} N(\theta)_l \big) \psi _k
\Big\}  
\Big] 
\\
& \qquad \qquad 
\cdot \nabla e^{(t-\tau)\Delta _D} h_j ~{\rm d}x {\rm d}\tau
\\
=: &
I_{\rm nonlinear}^1 + I_{\rm nonlinear}^2 ,
\end{split}
\]
where we denote by $I_{\rm nonlinear}^1, I_{\rm nonlinear}^2$  
the terms corresponding to the conditions $k \geq l$ and $k < l$, respectively. 

\vskip2mm 

We start by estimating $I_{\rm nonlinear}^2$. 
In the case when $j \geq k$, we proceed as follows. 
\[
\begin{split}
& \Big| 
\sum_{j\in \mathbb Z} 
\sum_{ k \leq j} \sum _{l \leq k}
\int_0^t 
\int_{\Omega} 
\Big[  
\Big\{ (\nabla^\perp \Lambda_D^{-1} \psi_k) N(\theta)_l
          + \big(\nabla ^\perp \Lambda _D^{-1} N(\theta)_l \big) \psi _k
\Big\}  
\Big] 
\cdot \nabla e^{(t-\tau)\Delta _D} h_j ~{\rm d}x {\rm d}\tau
\Big|
\\
\leq 
& C 
\sum_{j\in \mathbb Z} 
\sum_{ k \leq j} \sum _{l \leq k}
\int_0^t \| \psi_k \|_{L^{2p}} \| N(\theta)_l \|_{L^{2p}} 
  \cdot 2^j e^{-c(t-\tau)2^{2j}} \| h_j \|_{L^{2p}} 
  ~d\tau 
\\
\leq 
& C 
\sum_{j\in \mathbb Z} 
\Big(\sup _{\tau \in (0,\delta]}\sum_{ k \leq j} \sum _{l \leq k}
  2^{(1-\frac{1}{p})k}\| \psi_k \|_{L^{2}} 
  \cdot 2^{(1-\frac{1}{p})l} \| N(\theta)_l \|_{L^{2}} 
\Big) 
  \cdot 2^j \cdot 2^{-2j} \| h_j \|_{L^{2p}} 
  ~d\tau 
\\
\leq 
& C 
\| \psi \|_{L^\infty (0,\delta ; \dot B^0_{2,\infty})}
\| N(\theta) \|_{L^\infty (0,\delta ; \dot B^0_{2,\infty})} 
\| h \|_{\dot B^{1-\frac{2}{p}}_{p',1}} .
\\
\leq 
& C 
\| \psi \|_{L^\infty (0,\delta ; \dot B^{-1+\frac{2}{p}}_{p,\infty})}
\| N(\theta) \|_{L^\infty (0,\delta ; L^2 )} 
\| h \|_{\dot B^{1-\frac{2}{p}}_{p',1}} .
\end{split}
\]
In the case when $j < k$, 
we use equality \eqref{0823-2} and 
the H\"older inequality to obtain 
the $L^p, L^{p'}$, and $ L^\infty$ norms for $\psi_k , N(\theta)_l$, and $ h_j$ 
as follows. 
\[
\begin{split}
& \Big| 
\sum_{j\in \mathbb Z} 
\sum_{ k > j} \sum _{l \leq k}
\int_0^t 
\int_{\Omega} 
\Big[  
\Big\{ (\nabla^\perp \Lambda_D^{-1} \psi_k) N(\theta)_l
          + \big(\nabla ^\perp \Lambda _D^{-1} N(\theta)_l \big) \psi _k
\Big\}  
\Big] 
\cdot \nabla e^{(t-\tau)\Delta _D} h_j ~{\rm d}x {\rm d}\tau
\Big|
\\
\leq 
& C 
\sum_{j\in \mathbb Z} 
\sum_{ k > j} \sum _{l \leq k}
\int_0^t 2^{-k} \| \psi_k \|_{L^{p}} \| N(\theta)_l \|_{L^{p'}} 
  \cdot 2^{2j} e^{-c(t-\tau)2^{2j}} \| h_j \|_{L^{\infty}}  
  ~d\tau 
\\
\leq 
& C 
\sum_{j\in \mathbb Z} 
\sum_{ k > j} \sum _{l \leq k}
\sup_{\tau \in (0,\delta]}
2^{-k} \| \psi_k \|_{L^{p}} 2^{2(\frac{1}{2}-\frac{1}{p'})l} 
  \| N(\theta)_l \|_{L^{2}} 
  \cdot 2^{\frac{2}{p'}j} \| h_j \|_{L^{p'}}  
\\
\leq 
& C 
\| \psi \|_{L^\infty (0,\delta ; \dot B^{-1+\frac{2}{p}}_{p,\infty})}
\| N(\theta) \|_{L^\infty (0,\delta ; L^2 )} 
\| h \|_{\dot B^{1-\frac{2}{p}}_{p',1}} .
\end{split}
\]
Furthermore, we can  estimate $I_{\rm nonlinear}^2$ analogously,  
and thus, we obtain \eqref{0823-1}. 

\vskip3mm 

\noindent {\bf Step 4. } 
By the smallness in \eqref{0826-1} and the inequalities 
in Steps 2 and 3, it follows from the same argument as in 
\cite{IwUe-2024}, 
there exists $\delta > 0$ such that 
$\psi = 0 $ in $\dot B^0_{2,\infty}$ for all $t \in [0,\delta]$, 
which implies that $\psi = 0$ in $L^2$ for all $t \in [0,\delta]$.  
The uniqueness in the interval $[\delta ,T]$ can be proved 
by a contradiction argument. In fact, 
define
\begin{equation}\notag
 \tau^* =  \sup \big\{ \tau \in [0,T) \ \big| \ \lVert \theta(t, \cdot) - \tilde{\theta}(t, \cdot) \rVert_{L^2} = 0 \ \text{for} \ t \in [0, \tau] \big\}.
\end{equation}
If $\tau^* = T$, then the proof is completed.  
Assume that $\tau^* < T$. From the continuity in time of $\theta$ and $\tilde{\theta}$, it follows that $\theta(\tau^*) = \tilde{\theta}(\tau^*)$. 
Treating $\tau^*$ as the initial time, we argue as before, and 
the uniqueness holds in $[\tau^*, \tau^*+\delta']$ for some $\delta'>0$. This contradicts the definition of $\tau^*$.

 \vspace{5mm}

 \noindent
{\bf Data availability statement}. This manuscript has no associated data.

\noindent
{\bf Conflict of Interest. }
The author declares that he has no conflict of interest.

\end{document}